\documentclass [11pt,oneside]{article}

\reversemarginpar \textwidth=14cm \textheight=19cm

\usepackage{amssymb} \usepackage{amsfonts} \usepackage{amsmath}
\usepackage{amsthm} \usepackage{epsfig}
\usepackage{graphics}

\newtheorem{lemma}{Lemma}[section] 
\newtheorem{teo}[lemma]{Theorem}
\newtheorem{rem}[lemma]{Remark} 
\newtheorem{prop}[lemma]{Proposition}
\newtheorem{cor}[lemma]{Corollary}

\newcommand{\matR} {\ensuremath {\mathbb{R}}}

\newcommand{\matZ} {\ensuremath {\mathbb{Z}}}
\newcommand{\matC} {\ensuremath {\mathbb{C}}}

\newcommand{\matH} {\ensuremath {\mathbb{H}}}

\newcommand{\calI} {\ensuremath {\mathcal{I}}}

\newcommand{\calM} {\ensuremath {\mathcal{M}}}

\newcommand{\calP} {\ensuremath {\mathcal{P}}}
\newcommand{\calT} {\ensuremath {\mathcal{T}}}

\newcommand{\calO}{\ensuremath {\mathcal{O}}}

\newfont{\Got}{eufm10 scaled 1200}

\font\titsc=cmcsc10 scaled 1200

\newcommand{\finedimo}{{\hfill\hbox{$\square$}\vspace{2pt}}}

\author{Roberto \titsc{Frigerio}}

\title{An infinite family of hyperbolic graph complements in $S^3$}

\begin{document}

\maketitle

\begin{abstract}
For any $g\geqslant 2$ we construct a graph $\Gamma_g\subset
S^3$ whose exterior $M_g=S^3\setminus N(\Gamma_g)$
supports a complete finite-volume hyperbolic
structure with one toric cusp and a connected
geodesic boundary of genus $g$. We compute the canonical
decomposition and the isometry group of $M_g$, showing
in particular
that any self-homeomorphism of $M_g$ extends to a 
self-homeomorphism of the pair $(S^3,\Gamma_g)$, and that 
$\Gamma_g$ is chiral. 
Building on a result of Lackenby~\cite{lackenby}
we also show that any non-meridinal Dehn filling of $M_g$ 
is hyperbolic, thus getting an infinite family of 
graphs in $S^2\times S^1$ whose exteriors support a hyperbolic
structure with geodesic boundary. 
\vspace{4pt}

\end{abstract}

\section{Preliminaries and statements}

In this paper we introduce an infinite class $\{\Gamma_g,\ g
\geqslant 2\}$ of graphs in $S^3$
whose exteriors support a complete finite-volume hyperbolic structure
with geodesic boundary. Any $\Gamma_g$ has two connected components, one
of which is a knot.
We describe some geometric and
topological properties of the $\Gamma_g$'s, and we show that
for any $g\geqslant 2$
any non-meridinal Dehn-filling of the torus boundary of 
the exterior of $\Gamma_g$ gives a compact hyperbolic manifold
with geodesic boundary.

\paragraph{Definition of $\Gamma_g$ and hyperbolicity}
We say that a compact orientable $3$-manifold is \emph{hyperbolic}
if, after removing the boundary tori, we get a complete
finite-volume hyperbolic $3$-manifold with geodesic boundary.
Let $\Gamma$ be a graph in a closed $3$-manifold $M$ and let
$N(\Gamma)\subset M$ be an open regular neighbourhood
of $\Gamma$ in $M$. We say that $\Gamma$ is \emph{hyperbolic}
if $M\setminus N(\Gamma)$ is hyperbolic. If so, Mostow-Prasad's
Rigidity Theorem (see~\cite{FriPe,Fri} for a proof in the case with 
non-empty geodesic boundary) ensures that the 
complete finite-volume hyperbolic
structure with geodesic boundary 
on $M\setminus N(\Gamma)$ is unique up to isometry.

\begin{figure}
\begin{center}
\input{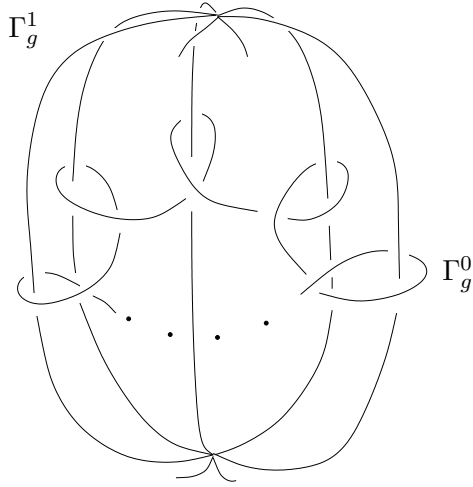}
\caption{\small{$\Gamma_g$ has two components: $\Gamma_g^0$ is a knot,
while $\Gamma_g^1$ is a graph with $g+1$ edges 
and two vertices.}}\label{graph:fig}
\end{center}
\end{figure}

\begin{figure}
\begin{center}
\includegraphics{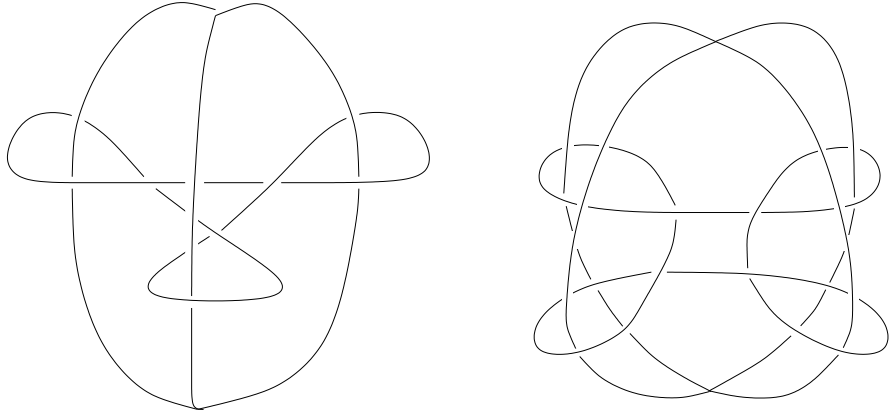}
\caption{\small{The graphs $\Gamma_2$ (on the left) and $\Gamma_3$
(on the right).}}\label{g2g3:fig}
\end{center}
\end{figure}

For any integer $g\geqslant 2$
let $\Gamma_g\subset S^3$ be the graph shown in Fig.~\ref{graph:fig}
(the graphs $\Gamma_2$ and $\Gamma_3$ are shown in Fig.~\ref{g2g3:fig}).
Let us denote by $\Gamma_g^0$ and $\Gamma_g^1$ the connected components
of $\Gamma$, where $\Gamma_g^0\cong S^1$ and
$\Gamma_g^1$ has two vertices and $g+1$ edges.
We also put $M_g=S^3\setminus N(\Gamma_g)$, 
$\partial_0 M_g=\partial N(\Gamma_g^0)$ and
$\partial_1 M_g=\partial N(\Gamma_g^1)$. 
Recall that if $M$ is a compact $3$-manifold
with $\partial M=\partial_0 M\sqcup \partial_1 M$,
the \emph{Heegaard genus} of $(M,\partial_0 M,\partial_1 M)$
is the minimal genus of a surface that splits $M$ as $C_0\sqcup C_1$,
where $C_i$ 
is obtained by attaching $1$-handles either
to the ``internal'' side of a collar of 
$\partial_i M$, if $\partial_i M\neq\emptyset$,
or to a $0$-handle if $\partial_i M=\emptyset$ 
(so $C_i$ is a handlebody in the latter case). 
In Section~\ref{triageo:sec}
we prove the following:

\begin{teo}\label{graph:teo} 
The graphs $\Gamma_g,\ g\geqslant 2$ are hyperbolic and the
Heegaard genus of $(M_g,\partial_0 M_g,\partial_1 M_g)$ is
$g+1$.
Moreover, the hyperbolic volume of $M_g$ grows linearly
with $g$ as follows:
\[
\lim_{g\to\infty} \frac{\mathrm{vol}(M_g)}{g}=5.419960359\ldots
\]
\end{teo}

\begin{figure}
\begin{center}
\includegraphics{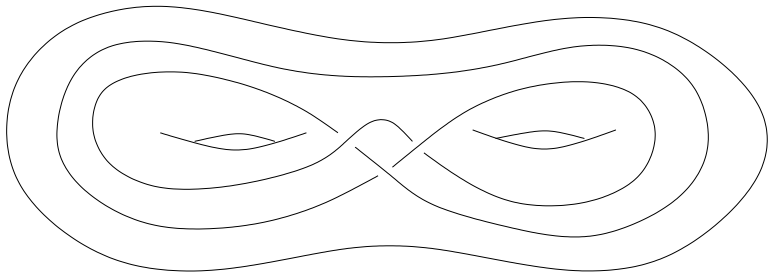}
\caption{\small{$M_2$ as a knot exterior in the handlebody 
of genus $2$.}}\label{adams:fig}
\end{center}
\end{figure}

\begin{rem}
\emph{Since $\Gamma_g^1$ is unknotted in $S^3$, the manifold
$M_g$ is the exterior of a knot in the handlebody of genus $g$.
The knot giving $M_2$, which is shown in 
Fig.~\ref{adams:fig}, was first introduced in~\cite{adams},
where it was proved to be hyperbolic by means of Thurston's
Hyperbolization Theorem for Haken manifolds.}
\end{rem}

\paragraph{Canonical decomposition and isometries}
In order to give $M_g$ a hyperbolic structure we 
geometrize a topological triangulation of $M_g$ by choosing
suitable shapes for the tetrahedra involved. More precisely,
let $\Delta$ denote the standard tetrahedron,
and let $\dot\Delta$ be $\Delta$ with its vertices removed.
An \emph{ideal triangulation} of a compact $3$-manifold $M$ 
with non-empty boundary is a realization of the interior
of $M$ as a gluing of a finite number of copies of
$\dot{\Delta}$, induced by a simplicial face-pairing
of the corresponding $\Delta$'s. In Section~\ref{triageo:sec}
we construct for any
$g\geqslant 2$ an ideal triangulation
$\calT_g$ of $M_g$ by $2g+2$ tetrahedra and we realize
the simplices of $\calT_g$ as geodesic polyhedra in $\matH^3$,
in such a way that the hyperbolic structure defined on them
extends to the whole of $M_g$.

Kojima proved in~\cite{Kojima} that
every complete finite-volume hyperbolic manifold
with non-empty geodesic boundary admits a \emph{canonical
decomposition} into geometric polyhedra.
The following result is proved in Section~\ref{canonical:sec}:

\begin{prop}\label{canonical:prop}
The canonical decomposition of $M_g$ is given by $\calT_g$.
\end{prop}

Let $\mathrm{Iso}(M_g)$ be the group of isometries of $M_g$,
let $\mathrm{Aut}(\calT_g)$ be the group of combinatorial 
automorphisms of $\calT_g$ and let $\calM (M_g)$ be the group
of homotopy classes of self-homeomorphisms of $M_g$.
An easy application of Mostow-Prasad's Rigidity Theorem gives
$\mathrm{Iso}(M_g)\cong\mathrm{Aut}(\calT_g)\cong\calM (M_g)$.

An oriented manifold $M$ is \emph{chiral} if it does not
admit an orientation-reversing self-homeomorphism.
Let $\Gamma\subset S^3$ be a graph and let $\Gamma'$ be the mirror
image of $\Gamma$. We say that $\Gamma$ is \emph{chiral} if
there does not exist an orientation-preserving  
homeomorphism between $(S^3,\Gamma)$ and $(S^3,\Gamma')$.
The definitions just given imply that $\Gamma$ is chiral
if $S^3\setminus\Gamma$ 
is chiral.
Using Proposition~\ref{canonical:prop} we will prove the 
following:

\begin{teo}\label{isometry:teo}
Let $D_n$ be the dihedral group of order $2n$, \emph{i.e.}~$D_n=
\langle r,s| \ r^n=s^2=1, rs=sr^{-1}\rangle$, and let
$g\geqslant 2$. Then:
\begin{enumerate}
\item
Any element of $\mathrm{Iso}(M_g)$
extends to a self-homeomorphisms
of $(S^3,\Gamma_g)$;
\item
The group $\mathrm{Iso}(M_g)$ is isomorphic to
$D_{g+1}$; 
\item
All the elements of $\mathrm{Iso}(M_g)$
are orientation-preserving; 
\item
The manifold $M_g$ and the graph $\Gamma_g$ are
chiral.
\end{enumerate}
\end{teo}

\paragraph{Dehn fillings}
Recall that a \emph{slope} on a torus is an isotopy class
of simple closed curves. 
For any $g\geqslant 2$, we denote by $s_g^m$ the \emph{meridinal}
slope on the torus boundary of $M_g$, \emph{i.e.}~the unique slope
on $\partial N(\Gamma_g^0)$ which bounds a disc in
$\overline{N(\Gamma_g^0)}$. For any slope $s$ in
$\partial N(\Gamma_g^0)$ we denote by
$M_g(s)$ the manifold obtained by Dehn-filling
the torus boundary of $M_g$ along $s$.
The next result is proved in Section~\ref{Dehn:sec}:

\begin{teo}\label{Dehn:teo}
Let $g\geqslant 2$ and let $s\neq s_g^m$ be a slope on the torus
boundary of $M_g$. Then $M_g(s)$ is hyperbolic and has Heegaard 
genus equal to $g+1$.
\end{teo}

\begin{figure}
\begin{center}
\input{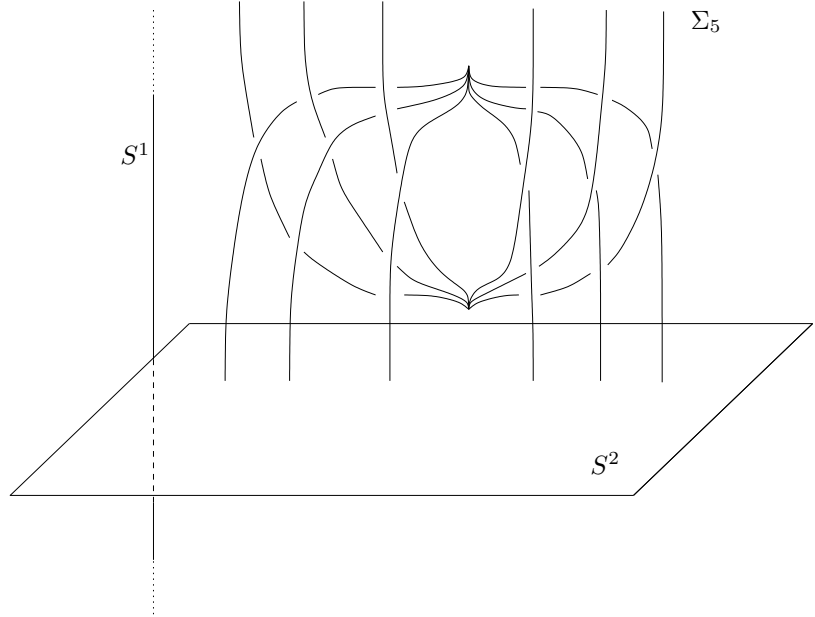}
\caption{\small{$\Sigma_g$ has $g+1$ edges and two vertices.
The picture shows the mirror image 
of $\Sigma_5$ as a graph in $\matR^2\times\matR\subset
S^2\times S^1$.}}\label{s2s1:fig}
\end{center}
\end{figure}

Since $\Gamma_g^0$ is unknotted in $S^3$, performing
$(0,1)$-Dehn surgery on the boundary torus of $M_g$
we get the exterior of a graph in $S^2\times S^1$. 
So Theorem~\ref{Dehn:teo} easily implies the following:

\begin{cor}\label{s2s1:cor}
Let $\Sigma_g$ be the graph in $S^2\times S^1$ shown
in Fig.~\ref{s2s1:fig}.
Then $\Sigma_g$ is a tunnel number one 
hyperbolic graph
for any $g\geqslant 2$.
\end{cor}
 
\begin{rem}
\emph{Theorem~\ref{Dehn:teo} also implies the following result,
which generalizes Corollary~\ref{s2s1:cor}: For any $g\geqslant 2$,
the $(p,q)$-Dehn surgery on $\Gamma_g^0$ yields a 
hyperbolic tunnel number one graph in the lens space $L_{p,q}$.}
\end{rem}

\section{Triangulations and hyperbolicity}\label{triageo:sec}

In this section we prove Theorem~\ref{graph:teo}.
To this aim we construct an ideal triangulation of $M_g$
and we prove that the tetrahedra of such triangulation
can be given hyperbolic structures which match under the gluings.

\begin{figure}
\begin{center}
\input{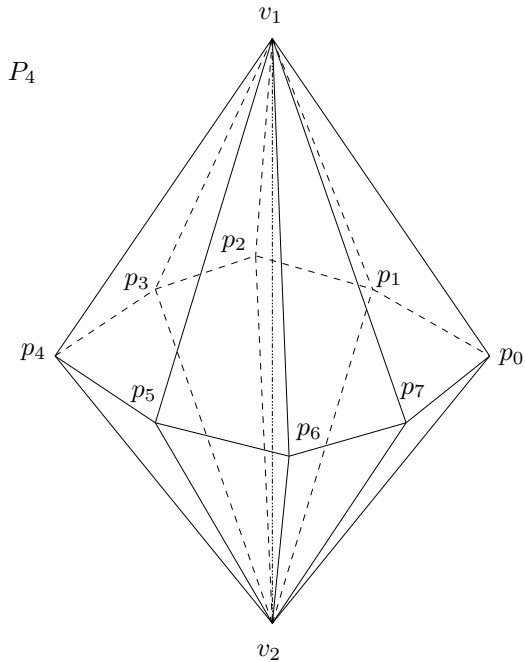}
\caption{\small{Gluing in pairs the faces of this
ideal double cone we get (the internal part of) 
$M_3$.}}\label{pir:fig}
\end{center}
\end{figure}

\paragraph{Defining the ideal triangulation}
For any $n\geqslant 3$, let $P_n$ be the solid double cone 
based on the regular $2n$-gon, and let $\dot{P}_n$ be 
$P_n$ with its vertices removed. We fix notation as suggested
in Fig.~\ref{pir:fig}, viewing $\mathrm{mod}\ 2n$
the  index $i$ of the $p_i$'s.
Let $X_n$ be the topological space obtained by gluing
the faces of $\dot{P}_n$ according to the following rules:
\begin{itemize}
\item
For any $i=0,2,\ldots,2n-4,2n-2$, the face $v_1 p_i p_{i+1}$
is identified with the face $p_{i+1} p_{i+2} v_2$
(with $v_1$ identified to $p_{i+1}$, $p_i$ to $p_{i+2}$ and
$p_{i+1}$ to $v_2$);
\item 
For any $i=1,3,\ldots,2n-3,2n-1$, the face $v_1 p_i p_{i+1}$
is identified with the face $p_{i+2} v_2 p_{i+1}$
(with $v_1$ identified to $p_{i+2}$, $p_i$ to $v_2$ and
$p_{i+1}$ to $p_{i+1}$).
\end{itemize}
The proof of the
following proposition will be sketched in the next paragraph.

\begin{prop}\label{triafond:prop}
For any $n\geqslant 3$,
$X_n$ is homeomorphic to the interior of $M_{n-1}$.
\end{prop}

Before discussing the proof of
Proposition~\ref{triafond:prop}, let us observe
that we can subdivide $P_{g+1}$ into $2g+2$ tetrahedra by adding
the ``vertical'' edge $v_1 v_2$. Such tetrahedra give 
an ideal triangulation of $M_{g}$, which we denote
from now on by $\calT_g$. The \emph{incidence number}
of an edge in a triangulation is the number of tetrahedra
incident to it (with multiplicity).
The definition of $\calT_g$ 
readily implies the following:

\begin{lemma}\label{edgeorder:lem}
For any $g\geqslant 2$, the triangulation $\calT_g$ has 
$g+3$ edges $e_0,\ldots,e_{g+2}$ such that:
\begin{itemize}
\item
for any $k=0,\ldots,g$, the edge $e_k$ is the projection in $M_g$ of
the edges
\[
v_1 p_{2k},\  p_{2k}p_{2k+1},\  p_{2k+1}p_{2k+2},\ p_{2k+2}
v_2\subset \dot{P}_{g+1},
\] 
and has incidence number $6$;  
\item
$e_{g+1}$ is the projection in $M_g$ of the edges 
$\{v_i p_j\subset \dot{P}_{g+1},\ i=1,2,\ j\ \mathrm{odd}\}$,
and has incidence number $4g+4$;
\item
$e_{g+2}$ is the projection in $M_g$ of $v_1 v_2\subset \dot{P}_{g+1}$,
and has incidence number $2g+2$.
\end{itemize}
\end{lemma}
\finedimo

\paragraph{Constructing the ideal triangulation}
In this paragraph we sketch the proof of 
Proposition~\ref{triafond:prop}. To this aim we apply 
(a slight generalization of)
the algorithm producing ideal triangulations for link complements
in $S^3$ described in~\cite{Pe}. Such an algorithm can be easily
modified in order to work with graphs rather than with links.
\begin{figure}
\begin{center}
\includegraphics{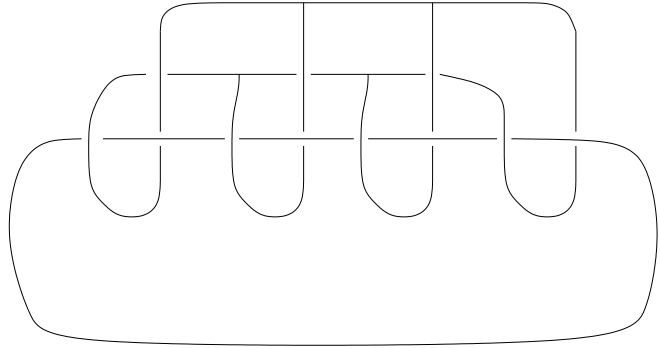}
\caption{\small{A projection of the graph $\Gamma'_3$.}}\label{proj:fig}
\end{center}
\end{figure}
Let us start with the following:
\begin{rem}\label{proj:rem}
\emph{Let $\Gamma'_g$ be the graph shown in Fig.~\ref{proj:fig}.
Then the complement of $\Gamma'_g$ is homeomorphic to the complement
of $\Gamma_g$, \emph{i.e.} $S^3\setminus N(\Gamma'_g)\cong M_g$.}
\end{rem}
Applying Petronio's algorithm to the projection
of $\Gamma'_g$ shown in Fig.~\ref{proj:fig} we obtain
the gluing diagrams shown in Fig.~\ref{diag:fig}.
Such diagrams encode the combinatorial rule which defines
the face-pairing of $P_{g+1}$ described in the previous
paragraph (see~\cite{Pe} for the details). This implies
Proposition~\ref{triafond:prop}.

\begin{figure}
\begin{center}
\input{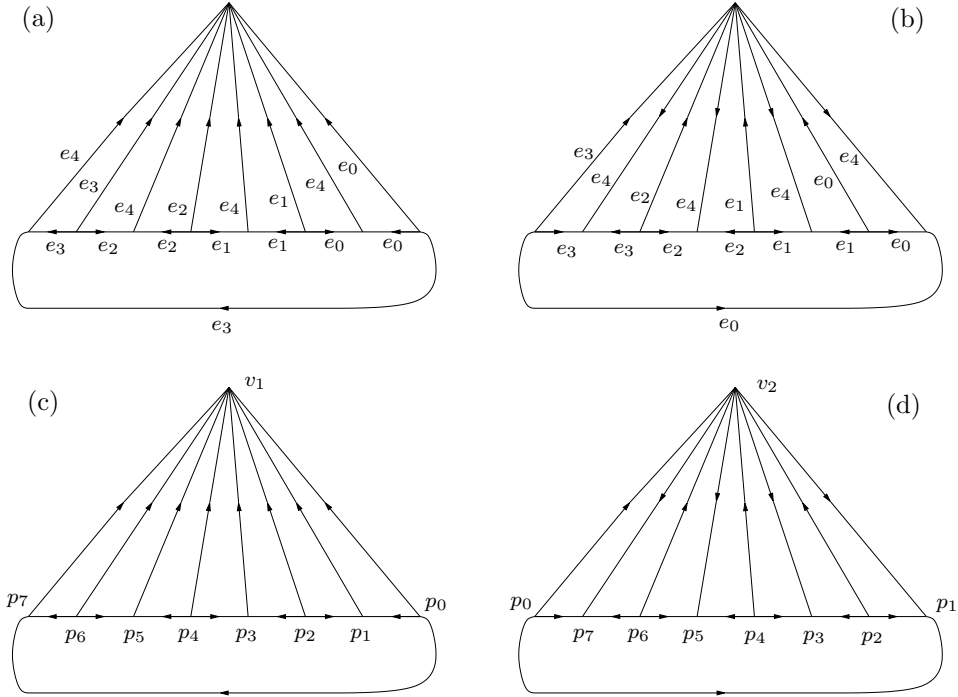}
\caption{\small{Applying Petronio's algorithm to the projection
in Fig.~\ref{proj:fig} we obtain the description of a cellularization
of $S^3\setminus \Gamma'_3$ which contains two
$3$-cells: in (a) and (c) it is described the
boundary of the upper $3$-cell, while in 
(b) and (d) we show the boundary of
the lower one. Labels for egdes and vertices agree
with notation in Fig.~\ref{pir:fig} and in 
Lemma~\ref{edgeorder:lem}.}}\label{diag:fig}
\end{center}
\end{figure}

\paragraph{Geometric tetrahedra}
In order to construct a hyperbolic structure on $M_g$
we will realize the tetrahedra of $\calT_g$ as 
geometric blocks in $\matH^3$. To describe the blocks to be used
we need some definitions.

A \emph{partially truncated tetrahedron} is a pair $(\Delta,\calI)$,
where $\Delta$ is a tetrahedron and $\calI$ is a set of vertices of $\Delta$,
which are called \emph{ideal vertices}. In the sequel we will always
refer to $\Delta$ itself as a partially truncated tetrahedron,
tacitly implying that $\calI$ is also fixed. The \emph{topological
realization} $\Delta^{\!\ast}$ of $\Delta$ is obtained by removing from $\Delta$ 
the ideal vertices and small open stars of the non-ideal ones.
We call \emph{lateral hexagon} and \emph{truncation triangle} the intersection
of $\Delta^{\!\ast}$ respectively with a face of $\Delta$ and with the link in
$\Delta$ of a non-ideal vertex. The edges of the truncation triangles,
which also belong to the lateral hexagons, are called \emph{boundary edges}, 
and the other edges of $\Delta^{\!\ast}$ are called \emph{internal edges}.
Note that, if $\Delta$ has ideal vertices, 
a lateral hexagon of $\Delta^{\!\ast}$
may not quite be a hexagon, because some of its (closed) edges may be missing.
A \emph{geometric realization} of $\Delta$ is an embedding of $\Delta^{\!\ast}$
in $\matH^3$ such that the truncation triangles are geodesic triangles,
the lateral hexagons are geodesic polygons with ideal vertices corresponding
to missing edges, and the truncation triangles and lateral hexagons 
lie at right angles to each other.

\paragraph{Consistency}
For any $g\geqslant 2$ let us set
\[
\alpha_g=\pi/(2g+2),\ \beta_g=2\alpha_g,
\ \gamma_g=\arccos\left( (2\cos \alpha_g)^{-1}\right),
\ \delta_g=\pi-2\gamma_g.
\]
For any $i=0,\ldots,2g+1$ let
$\Delta_i$ be the tetrahedron in $P_{g+1}$ 
with vertices $v_1, v_2, p_i, p_{i+1}$.
We realize the simplices of 
$\calT_g$ as partially truncated tetrahedra 
\begin{figure}[h]
\begin{center}
\input{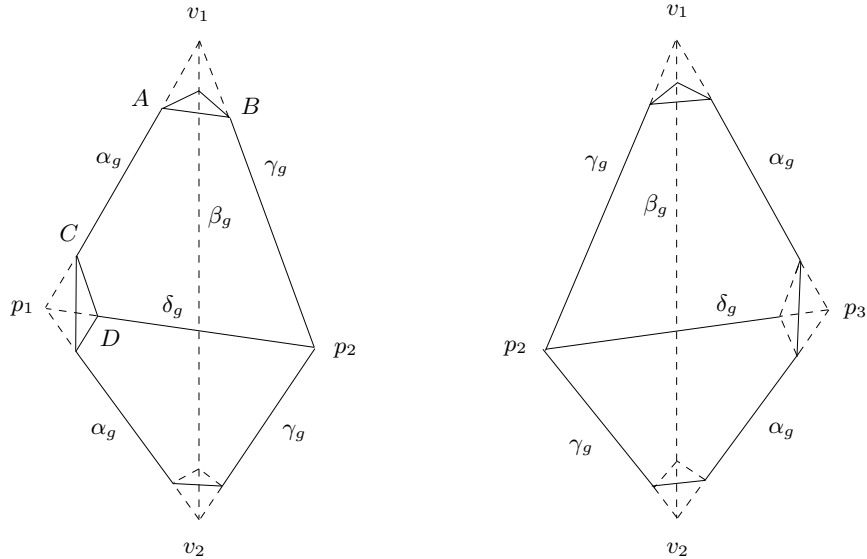}
\caption{\small{The dihedral angles along the edges
of the tetrahedra in $\calT_g$: to ensure that the matching faces
can be glued by isometries we impose that $AB$ and $CD$ have
the same length.}}\label{tetra:fig}
\end{center}
\end{figure}
as follows
(see Fig.~\ref{tetra:fig}):
\begin{itemize}
\item
For any $i=0,2,\ldots,2g-2,2g$, we declare $p_i$ 
to be the only ideal vertex of $\Delta_i$ and
we set $\Delta_i^{\!\ast}$ to be
the geometric realization of $\Delta_i$ with dihedral angles
$\delta_g$ along $p_i p_{i+1}$, $\beta_g$ along $v_1 v_2$,
$\gamma_g$ along $p_{i} v_1$ and $p_{i} v_2$, and $\alpha_g$
along $p_{i+1} v_1$ and $p_{i+1} v_2$;
\item
For any $i=1,3,\ldots,2g-1,2g+1$, we declare $p_{i+1}$ to be
the only ideal vertex of $\Delta_i$ and we set $\Delta_i^{\!\ast}$ to be
the geometric realization of $\Delta_i$ with dihedral angles
$\delta_g$ along $p_i p_{i+1}$, $\beta_g$ along $v_1 v_2$,
$\alpha_g$ along $p_i v_1$ and $p_i v_2$, and $\gamma_g$
along $p_{i+1} v_1$ and $p_{i+1} v_2$.
\end{itemize}
Existence and uniqueness 
of such geometric realizations are proved in~\cite{FriPe},
where
it is also shown that the hyperbolic structure given
on the tetrahedra of $\calT_g$ extends to the whole of $M_g$ if and only
if the matching boundary edges have the same length and the total dihedral 
angle around each internal edge is $2\pi$. Our choice of 
angles is such that all the conditions on dihedral angles and several 
conditions on boundary lengths are trivially satisfied.
The only non-trivial condition to be imposed in order
to ensure geometricity of $\calT_g$ is that
the edges $AB$ and $CD$ in Fig.~\ref{tetra:fig} should have the
same length. This requirement
translates into the following equation:
\begin{equation}\label{consistency:eq}
\frac{\cos \gamma_g \cos \alpha_g+\cos \beta_g}
{\sin \gamma_g \sin \alpha_g}
=\frac{\cos \delta_g \cos \alpha_g+\cos \alpha_g}
{\sin \delta_g \sin \alpha_g},
\end{equation}
which is solved by the choosen values for
$\alpha_g,\beta_g,\gamma_g$ and $\delta_g$. We have thus
proved that the geometric realization of $\calT_g$
just described defines a (possibly incomplete) hyperbolic
structure on the whole of $M_g$.

\begin{figure}
\begin{center}
\input{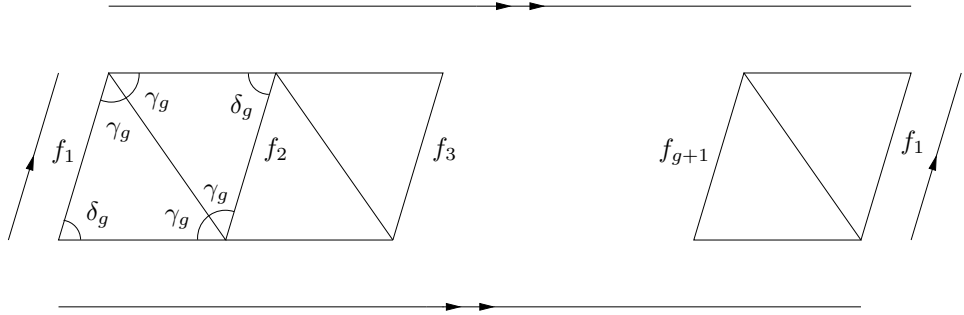}
\caption{\small{$\calT_g$ induces a tiling of the boundary
torus of $M_g$ by $2g+2$ isometric Euclidean isosceles
triangles.}}\label{tiling:fig}
\end{center}
\end{figure}

\paragraph{Completeness}
To check completeness of the hyperbolic structure defined in
the last paragraph
we have to determine the similarity structure it induces on the
boundary torus of $M_g$. By construction, the torus in $\partial M_g$ is
tiled by $2g+2$ Euclidean triangles as in Fig.~\ref{tiling:fig}.
This shows that the structure on the boundary torus is indeed 
Euclidean, so the hyperbolic structure constructed in the previous 
paragraph is complete.

\paragraph{Volume and Heegaard genus}
The notion of partially truncated tetrahedron introduced
before admits a generalization~\cite{FriPe} 
to the case in which some internal
edge degenerates into an ideal point, becoming a so-called
\emph{length-$0$ edge}. 
If the dihedral angle along a length-$0$ edge is declared
to be equal to $0$, then the volume
is a continuous function of the dihedral angles of \emph{generalized}
partially truncated tetrahedra (see~\cite{Ushijima}).

Let now $w_1,\ldots,w_4$ be the vertices of a tetrahedron $\Delta$, and suppose
that $w_1$ is the unique ideal vertex of $\Delta$. Let $\Delta_g^{\!\ast}$
be the geometric realization of $\Delta$ parametrized by the following
dihedral angles:
\begin{equation}\label{angles:eq}
\angle w_1 w_2=\angle w_1 w_3=\gamma_g,\ \angle w_1 w_4=\delta_g,
\ \angle w_2 w_4=\angle w_3 w_4=\alpha_g,\ \angle w_2 w_3=\beta_g,
\end{equation}
and set $V_g=\mathrm{Vol}(\Delta_g^{\!\ast})$. Let also $\Delta_\infty^{\!\ast}$
be the geometric (generalized) realization of $\Delta$ with dihedral
angles equal to 
$\pi/3$ along the edges emanating from $w_1$ and
equal to $0$ along the other three edges, and set
$V_\infty=\mathrm{Vol}(\Delta^{\!\ast}_\infty)$. Continuity of
volume of generalized partially truncated tetrahedra as a function of
their dihedral angles implies that
\[
\lim_{g\to\infty} \frac{\mathrm{Vol}(M_g)}{g}=
\lim_{g\to\infty} \frac{(2g+2)V_g}{g}=
2\lim_{g\to\infty} V_g=2 V_\infty=5.419960359\ldots
\]

Finally, the genus of $(M_g,\partial_0 M_g,\partial_1 M_g)$ is
of course at least $g$, and it is actually at most $g+1$
because the boundary of a regular neighbourhood of
$\partial_1 M_g\cup e_{g+1}$ is easily seen to be a Heegaard
surface. If the Heegaard genus of $M_g$ were $g$, then
 $\partial_1 M_g$ should be compressible in $M_g$, against
the hyperbolicity of $M_g$.

\section{Canonical decomposition}\label{canonical:sec}
Kojima proved in~\cite{Kojima} that a complete finite-volume
hyperbolic manifold $M$ with non-empty geodesic boundary
admits a \emph{canonical decomposition} into partially truncated 
polyhedra (an obvious generalization of a partially truncated 
tetrahedron). This decomposition is obtained by projecting first 
to $\matH^3$ and then to $M$ the faces of the convex hull of a certain
family $\calP$ of points in Minkowsky $4$-space. This family $\calP$
splits as $\calP'\sqcup\calP''$, with $\calP'$ consisting
of the points on the hyperboloid $||x||^2=+1$ which are dual to
the hyperplanes giving $\partial \widetilde{M}$, where 
$\widetilde{M}\subset\matH^3$ is a universal cover of $M$.
The points in $\calP''$ lie in the light-cone, and they are the duals
of horoballs projecting in $M$ to Margulis neighbourhoods of the cusps.
The choice of these Margulis neighbourhoods is somewhat tricky,
and carefully explained in~\cite{FriPe}. 
It will be sufficient for our present purposes to know
that any choice of sufficiently small Margulis neighbourhoods leads to a
set $\calP''$ which works. 
In the sequel we will denote by $\calO$ the union of 
sufficiently small Margulis neighbourhoods of the cusps.

\paragraph{Tilts}
Suppose a decomposition $\calT$ of 
$M$ by partially truncated tetrahedra is given. 
The matter of deciding if $\calT$ is the canonical Kojima decomposition
of $M$ is faced using the \emph{tilt formula}~\cite{weeks:tilt, Ush, FriPe},
that we now briefly describe.

Let $\sigma$ be a $d$-simplex in $\calT$ and $\tilde{\sigma}$ be
a lifting of $\sigma$ to $\widetilde{M}\subset\matH^3$. 
To each end of $\tilde{\sigma}$ 
there corresponds (depending on the nature of the end) 
one horoball in the lifting of $\calO$ 
or one component of the geodesic boundary of $\widetilde{M}$, 
so $\tilde{\sigma}$
determines $d+1$ points
of $\calP$. Now let two 
tetrahedra $\Delta^{\!\ast}_1$ and $\Delta^{\!\ast}_2$ 
share a $2$-face
$F^\ast$, and let 
$\widetilde{\Delta}^{\!\ast}_1,\widetilde{\Delta}^{\!\ast}_2$ 
and $\widetilde{F}^{\ast}$ be 
liftings of $\Delta^{\!\ast}_1,\Delta^{\!\ast}_2$ 
and $F^\ast$ to $\widetilde{M}\subset\matH^3$ 
such that $\widetilde{\Delta}^{\!\ast}_1
\cap\widetilde{\Delta}^{\!\ast}_2=\widetilde{F}^{\ast}$.
Let $\overline{F}$ be
the $2$-subspace in Minkowsky 4-space that
contains the three points of $\calP$ 
determined by $\widetilde{F}^{\ast}$. For $i=1,2$ let 
$\overline{\Delta}^{(F)}_i$ be the half-$3$-subspace bounded by 
$\overline{F}$ and containing the fourth point of $\calP$ determined
by $\widetilde{\Delta}^{\!\ast}_i$. 
Then one can show that $\calT$ is canonical if and only if, 
whatever $F^\ast\!\!,\Delta^{\!\ast}_1,\Delta^{\!\ast}_2$, the following holds:
\begin{itemize}
\item The half-$3$-subspaces
$\overline{\Delta}^{(F)}_1$ and $\overline{\Delta}^{(F)}_2$
lie on  distinct $3$-subspaces and their convex hull 
does not contain the origin of Minkowsky $4$-space.
\end{itemize}
The tilt formula computes 
a real number $t(\Delta^{\!\ast},F^\ast)$ describing the
``slope'' of $\overline{\Delta}^{(F)}$ in terms of 
the intrinsic geometry of $\Delta^{\!\ast}$ and $\calO$. More precisely,
one can translate the condition just stated 
into the inequality $t(\Delta^{\!\ast}_1,F^\ast)
+t(\Delta^{\!\ast}_2,F^\ast)< 0$.

\paragraph{The canonical decomposition of $M_g$}
Coming to the manifolds we are interested in, let $g\geqslant 2$,
let $\calT_g$ be the geometric triangulation
of $M_g$ we have described in the previous section 
and let $\calO$ be a suitable neighbourhood
of the cusp of $M_g$. It was shown in~\cite{FriPe} that $\calO$
determines a real number $r_\Delta(v)>0$ for the ideal vertex $v$
of any tetrahedron $\Delta$ in $\calT_g$.
This number $r_\Delta(v)$ represents the ``height'' 
of the trace in $\Delta$ near $v$ of $\partial\calO$
(except that $r_\Delta(v)\ll 1$ means that $\partial\calO$ is ``very'' high).
Looking at the intersection of $\calO$ with the tetrahedra
of $\calT_g$ it is easily seen that
$r_\Delta(v)$ has a certain fixed
value $r$ whenever $v$ is the ideal vertex of any $\Delta$ in $\calT_g$.

Let now $w_1,\ldots,w_4$ be the vertices of a 
geometric partially truncated tetrahedron $\Delta_g^{\!\ast}$
in $\calT_g$ and suppose
that $w_1$ is the unique ideal vertex of $\Delta_g^{\!\ast}$. 
Let the dihedral angles of $\Delta_g^{\!\ast}$ 
be as prescribed in equation~(\ref{angles:eq})
and let $r=r(w_1)$ be the parameter associated with the intersection
of $\calO$ with $\Delta_g^{\!\ast}$. For any $i=1,\ldots,4$ let
also $F^\ast_i$ be the face of $\Delta_g^{\!\ast}$ opposite to $w_i$.
Using the formulae given in~\cite{FriPe} 
we can compute the tilts of the geometric blocks
of $\calT_g$. Let us set
\begin{eqnarray*}
d_1&=&2r\cdot(\sin^2 \gamma_g\sin\delta_g)/
(2\sin\gamma_g\cos\alpha_g+\sin\delta_g\cos\beta_g)>0,\\
d_2=d_3&=&\cos^2 \alpha_g+\cos^2 \beta_g+\cos^2\gamma_g 
+ 2\cos \alpha_g\cos \beta_g\cos\gamma_g-1>0,\\
d_4&=&2\cos^2 \alpha_g+\cos^2 \delta_g+2\cos^2 
\alpha_g\cos\delta_g-1>0.
\end{eqnarray*}
Then there exists a constant $k_g>0$, depending only on $g$, such that  
\begin{eqnarray*}
t(\Delta_g^{\!\ast},F^\ast_1)&=&
d_1-
k_g\cdot\left(2\sqrt{d_2}\cos
\alpha_g+
\sqrt{d_4}\cos\beta_g\right),\\
t(\Delta_g^{\!\ast},F^\ast_2)=
t(\Delta_g^{\!\ast},F^\ast_3)&=&-d_1\cos\alpha_g+
k_g\cdot\left(\sqrt{d_2}
(1-\cos\delta_g)-\sqrt{d_4}\cos\gamma_g\right),\\
t(\Delta_g^{\!\ast},F^\ast_2)&=&-d_1\cos \beta_g+
k_g\cdot\left(-2\sqrt{d_2}\cos\gamma_g+\sqrt{d_4}\right).
\end{eqnarray*}
Now an easy computation shows that
\[
\sqrt{d_2}(1-\cos\delta_g)
-\sqrt{d_4}\cos\gamma_g=
-2\sqrt{d_2}\cos\gamma_g+
\sqrt{d_4}=0.
\]
This implies that if $r$ is small enough, then all the tilts
are negative, so $\calT_g$ is the Kojima decomposition of $M_g$
for any $g\geqslant 2$, and Proposition~\ref{canonical:prop} is proved.

\paragraph{Isometry group}
By Proposition~\ref{canonical:prop}, the isometry group of
$M_g$ is canonically isomorphic to the group $\mathrm{Aut}
(\calT_g)$ of combinatorial automorphisms of $\calT_g$. 
For any $g\geqslant 2$, let $r_g$ and $s_g$ be the unique combinatorial 
automorphisms of the double cone
$P_{g+1}$ such that:
\begin{eqnarray*}
r_g(p_i)=p_{i+2}\ \forall 
i\in\{0,\ldots,2g+1\},\ r_g(v_1)=v_1,\ r_g(v_2)=v_2,\\
s_g(p_i)=p_{-i} 
\ \forall i\in\{0,\ldots,2g+1\},\ s_g(v_1)=v_2,
\ s_g(v_2)=v_1.
\end{eqnarray*}
If we realize $P_{g+1}$ as a Euclidean
regular $(2g+2)$-gonal double cone, then 
$r_g$ is a rotation of order $g+1$
around the line containing $v_1$ and $v_2$, 
while $s_g$ is a rotation of order
$2$ having as axis the line through
$p_0$ and $p_{g+1}$. This easily implies that the following relations
hold:
\[
r_g^{g+1}=s_g^{2}=1,\ \ r_g\circ s_g=s_g\circ r_g^{-1}.
\]
Theorem~\ref{isometry:teo} is now readily deduced from the following:

\begin{prop}\label{iso:prop}
Both $r_g$ and $s_g$ induce combinatorial automorphisms of
$\calT_g$. Moreover, the group $\mathrm{Aut}
(\calT_g)$ is generated by $r_g$ and $s_g$ for any $g\geqslant 2$.
\end{prop} 
\noindent\emph{Proof:}
The first statement readily follows by a direct
computation. In order to prove the second statement,
for
any $g\geqslant 2$ let $H_g$ be the subgroup of
 $\mathrm{Aut}
(\calT_g)$ generated by $r_g$ and $s_g$.
Noting that $H_g$ acts transitively on the set
of tetrahedra of $\calT_g$, to conclude that
$\mathrm{Aut}(\calT_g)=H_g$ it is enough to show that
the stabilizer in $\mathrm{Aut}(\calT_g)$ of one
tetrahedron of $\calT_g$ is trivial. So
let $\Delta_0\subset P_{g+1}$ be the tetrahedron in $\calT_g$ 
with vertices $p_0,p_1,v_1,v_2$. 
We observe that $p_0$ and $p_1$ can be characterized 
as the only vertices of $\Delta_0$ satisfying the following
properties:
\begin{itemize}
\item
$p_0$ is asymptotic in $M_g$ to the boundary torus
(\emph{i.e.} it is ideal in the geometric realization);
\item
there exist edges $e,e'$ of $\Delta_0$ such that
$p_1=e\cap e'$ and $e,e'$ are projected in $M_g$
to the only edge of $\calT_g$ having order $4g+4$.
\end{itemize}
Let now $\varphi$
be an element in $\mathrm{Aut}(\calT_g)$
such that $\varphi(\Delta_0)=\Delta_0$.
The intrinsic description of $p_0$ and $p_1$ 
given above implies
that $\varphi(p_i)=p_i$ for $i=0,1$.
Now it is easily seen that the simplicial automorphism
of $\Delta_0$ which fixes $p_0,p_1$ and interchanges
$v_1$ with $v_2$ does not extend to an automorphism
of $\calT_g$, so $\varphi|_{\Delta_0}$ must be the
identity. This implies that
$\varphi=1$ in $\mathrm{Aut}(\calT_g)$, so the stabilizer
of $\Delta_0$ in $\mathrm{Aut}(\calT_g)$ is trivial,
and we are done.
\finedimo


\section{Dehn filling}\label{Dehn:sec}

This section is devoted to the proof of Theorem~\ref{Dehn:teo}.
A detailed analysis of the geometry of the tetrahedra of $\calT_g$ will
lead us to
an estimate on the size of the cusp of $M_g$. We will then apply
a theorem
of Lackenby~\cite{lackenby} to prove the desired result.

\paragraph{Maximal cusps, slopes and hyperbolicity}
Let $M$ be an orientable hyperbolic $3$-manifold and 
let $\dot{M}$ be $M$ with its boundary tori removed,
so that $\dot{M}$ admits by definition a complete finite-volume
hyperbolic structure with geodesic boundary and cusps.
We say that $\calO\subset \dot{M}$
is a \emph{horocusp section} for $M$ if 
$\dot{M}\setminus \calO$ is compact
and the preimage of $\calO$ in the universal 
covering of $\dot{M}$ is the disjoint union
of open horoballs. 
If $\calO$ is a horocusp section, then 
$\partial\calO$ is a union of (possibly touching or self-touching) 
tori and the 
hyperbolic structure of $\dot{M}$ induces
on $\partial\calO$ a well-defined
Euclidean metric.
A horocusp section for $M$ is 
\emph{maximal} if it is maximal among the horocusp
sections for $M$
with respect to inclusion. We observe that 
a maximal horocusp section always exists, and is unique
if $M$ has only one boundary torus.

Let now $T_1,\ldots,T_k$ be the boundary
tori of $M$, let $h\leqslant k$ and let
$s_i$ be a slope on $T_i$ for any $i=1,\ldots,h$.
We denote by $M(s_1,\ldots,s_h)$ the manifold
obtained by Dehn-filling $M$ along $s_1,\ldots,s_h$. 
If $\calO$ is a fixed horocusp section for $M$, then any
$s_i$ determines a well-defined isotopy class of Euclidean
geodesics on the corresponding component 
of $\partial\calO$. We denote by $L_\calO(s_i)$
the Euclidean length of such geodesics. The following theorem is
proved in~\cite{lackenby}:

\begin{teo}\label{lack:teo}
Let $M$ be hyperbolic with boundary
$\partial M=T_1\sqcup\ldots\sqcup T_k$ given by $k$ tori.
Let $h\leqslant k$ and $s_i$ be a slope on $T_i$
for any $i=1,\ldots,h$. Fix a horocusp
section $\calO$ for $M$ and suppose that $L_\calO(s_i)>6$
for any $i=1,\ldots,h$. Then $M(s_1,
\ldots,s_h)$ is irreducible, atoroidal and has infinite
word-hyperbolic fundamental group.
\end{teo}

\paragraph{Hyperbolic Dehn fillings}
Building on Theorem~\ref{lack:teo} we prove here the following:

\begin{prop}\label{lack:prop}
Let $M$ be hyperbolic with boundary
$\partial M=\Sigma_1\sqcup\ldots\sqcup\Sigma_r\sqcup
T_1\sqcup\ldots\sqcup T_k$ given by $r>0$ surfaces
of negative Euler characteristic and $k$ tori.
Let $h\leqslant k$ and $s_i$ be a slope on $T_i$
for any $i=1,\ldots,h$. Fix a horocusp
section $\calO$ for $M$ and suppose that $L_\calO(s_i)>6$
for any $i=1,\ldots,h$. Then $M(s_1,
\ldots,s_h)$ is hyperbolic.
\end{prop}
\noindent\emph{Proof:}
Let $DM$ be the double of $M$ along 
$\Sigma=\Sigma_1\sqcup\ldots\sqcup\Sigma_r$, \emph{i.e.} the
manifold obtained by mirroring $M$ along $\Sigma$.
Then $M$ is canonically embedded in $DM$, and
$DM$ admits an involution $\sigma$ which fixes
$\Sigma\subset M\subset DM$ 
and interchanges $M$ with its mirror copy.
Of course we have
$\partial DM=T_1\sqcup\ldots \sqcup T_k\sqcup \overline{T}_1\sqcup
\ldots\sqcup \overline{T}_k$, where $\overline{T}_i=\sigma(T_i)$ 
for any $i=1,\ldots,k$. Note also
that for any $i=1,\ldots,h$
the slope $s_i$ on $T_i$ determines a mirror slope
$\overline{s}_i=\sigma(s_i)$ on $\overline{T}_i$.
The manifold $DM$ can be given a  hyperbolic structure simply
by doubling the hyperbolic structure with geodesic boundary
of $M$. Then $\calO$ determines a horocusp
section $D\calO$ for $DM$ such that 
$L_{D\calO}(\overline{s}_i)=
L_{D\calO}(s_i)=
L_{\calO}(s_i)>6$ for any $i=1,\ldots,h$.
So Theorem~\ref{lack:teo} implies that
$DM(s,\overline{s}):=DM(s_1,\ldots,s_h,\overline{s}_1,
\ldots,\overline{s}_h)$ is irreducible
and atoroidal
and has infinite, word-hyperbolic fundamental group.

Let now $W$ be a compression disc
for $\Sigma$  
in $DM(s,\overline{s})$, so $\partial W=\Sigma\cap W$ 
does not bound a disc on $\Sigma$.
By possibly replacing $W$ with $\sigma(W)$ we can
suppose that $W$ is properly embedded in
$M(s_1,\ldots,s_h)$. If $W$ were non-separating in
$M(s_1,\ldots,s_h)$ then its double $DW$ would be a
non-separating sphere in $DM(s,\overline{s})$,
a contradiction since $DM(s,\overline{s})$
is irreducible.
So let $N_1$ and $N_2$ be the manifolds obtained by cutting
$M(s_1,\ldots,s_h)$ along $W$ and for $r=1,2$ let $DN_r$
be the double of $N_r$ along $\Sigma\cap N_r\subset \partial N_r$.
By construction the sphere $DW$ decomposes $DM(s,\overline{s})$
as the connected sum of $DN_1$ and $DN_2$, so $DN_{\overline{r}}$
is a $3$-ball for $\overline{r}=1$ or $2$. This easily implies
that $N_{\overline{r}}$ is $3$-ball and $\Sigma\cap N_{\overline{r}}$
is a disc on $\Sigma$ with boundary $\partial W$, a contradiction
since $W$ is a compression disc for $\Sigma$.
We have thus shown that $\Sigma$
is incompressible in $DM(s,\overline{s})$, so 
Thurston's Hyperbolization Theorem for Haken manifolds~\cite{thu2}
implies that $DM(s,\overline{s})$ is hyperbolic.

Of course $\sigma$ extends to an involution $\hat{\sigma}$ of 
$DM(s,\overline{s})$, which is in turn homotopic to
an involutive isometry $\overline{\sigma}$ by Mostow-Prasad's
Rigidity Theorem. A result of Tollefson~\cite{Tollefson}
now ensures that $\overline{\sigma}$ fixes a surface
$\Sigma'_1\sqcup\ldots\sqcup\Sigma'_r$ isotopic to
$\Sigma_1\sqcup\ldots\sqcup\Sigma_r$. Cutting 
$DM(s,\overline{s})$ along 
$\Sigma'_1\sqcup\ldots\sqcup\Sigma'_r$ we 
obtain isometric complete finite-volume 
hyperbolic manifolds with geodesic boundary
$Y$ and $\overline{Y}$. Since 
$M(s_1,\ldots,s_h)$ is homeomorphic to $Y$,
the conclusion follows.
\finedimo 

\begin{figure}[ht]
\begin{center}
\input{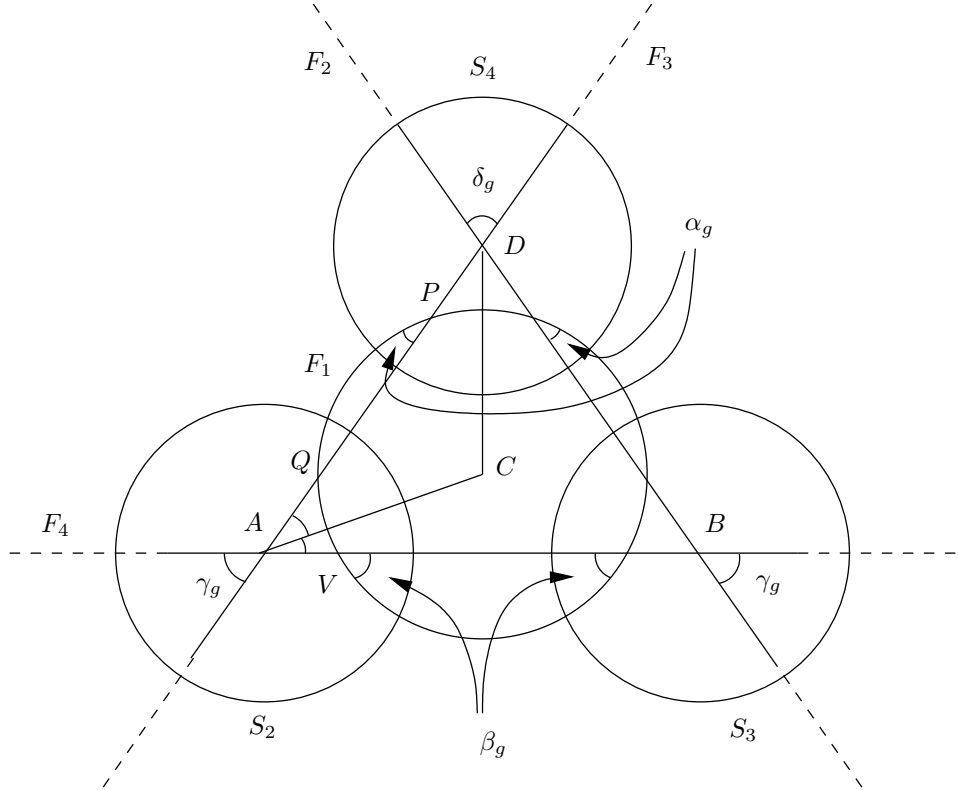}
\caption{\small{Notation for the proof of 
Lemma~\ref{primo:lem}.}}\label{pattern:fig}
\end{center}
\end{figure}

\paragraph{The maximal cusp of $M_g$}
Coming back to the case we are interested in,
we now want to determine the size of the maximal horocusp section
$\calO_g$ for $M_g$, $g\geqslant 2$. 
So let us fix $g\geqslant 2$ and
let $\Delta_g^{\!\ast}$ be a geometric partially truncated tetrahedron
in $\calT_g$ with vertices $w_1,w_2,w_3$ and $w_4$
and with dihedral angles as prescribed in 
equation~(\ref{angles:eq}).
We can realize $\Delta_g^{\!\ast}$ in the upper half-space model of
$\matH^3$ in such a way that $w_1$ is identified with $\infty$, 
so the truncation triangles corresponding to $w_2,w_3$ and $w_4$
lie on the hyperbolic planes bounded by $S_2,S_3$ and $S_4$ respectively, 
where $S_i$ is a Euclidean circle in
$\matC\times\{0\}\subset\partial\matH^3$ for $i=2,3,4$
(see Fig.~\ref{pattern:fig}).
We also denote by $F_1,F_2,F_3,F_4$ the traces at infinity
of the hyperbolic planes containing the faces of $\Delta_g^{\!\ast}$ with
vertices $w_2w_3w_4$, $w_1w_3w_4$, $w_1w_2w_4$, $w_1w_2w_3$ respectively.
By construction $F_1$ is a Euclidean circle perpendicular to 
$S_i$ for $i=2,3,4$,
while $F_j$ is a Euclidean line perpendicular to $S_l$ for $j,l\in\{2,3,4\}$,
$l\neq j$.
We denote by $C$ and $R$ the center and the radius of $F_1$,
and we also set
\[
A=F_3\cap F_4,\ B=F_2\cap F_4,\ D=F_2\cap F_3,
\ \{P,Q\}=F_1\cap F_3,\ V\in F_1\cap F_4.
\]
Let us observe that equation~(\ref{consistency:eq}) in 
Section~\ref{triageo:sec} implies that 
the Euclidean radii
of $S_2,S_3$ and $S_4$ are equal to the same
value $R'$. 
 
\begin{lemma}\label{primo:lem}
The following equalities hold:
\begin{eqnarray}
{R'}^2&=&R^2\cdot\frac{\cos^2\alpha_g+\cos^2\beta_g
+\cos^2\gamma_g+2\cos\alpha_g
\cos\beta_g\cos\gamma_g-1}{\sin^2\gamma_g},\label{first:eq}\\
\overline{AB}&=&2R'\cdot\frac{\cos\alpha_g+\cos\beta_g
\cos\gamma_g}{\sqrt{\cos^2\alpha_g+\cos^2\beta_g
+\cos^2\gamma_g+2\cos\alpha_g
\cos\beta_g\cos\gamma_g-1}}.\label{sec:eq}
\end{eqnarray}
\end{lemma} 
\noindent\emph{Proof:}
Let $\gamma'_g=\angle CAV,\ \gamma''_g=\angle CAQ$.
Applying the Euclidean Sine Rule to the triangles
$ACV$ and $ACQ$ we get
$\overline{AC}/\sin(\pi/2+\beta_g)=R/\sin \gamma'_g$ and
$\overline{AC}/\sin(\pi/2+\alpha_g)=R/\sin \gamma''_g$. So
we have $\sin\gamma''_g\cos\beta_g=\sin\gamma'_g\cos\alpha_g$.
Substituting the equality $\gamma''_g=\gamma_g-\gamma'_g$
in this equation and dividing by $\sin\gamma'_g$ we obtain:
\begin{equation}\label{cotang:eq}
\cot \gamma'_g=\frac{\cos\alpha_g+\cos\beta_g\cos\gamma_g}
{\cos\beta_g\sin\gamma_g}.
\end{equation}
Combining this relation with the equalities
\[
\sin\gamma'_g=(1+\cot^2 \gamma'_g)^{-1/2},
\ \overline{AC}
=R\cdot\cos\beta_g/\sin\gamma'_g,\ {R'}^2=\overline{AC}^2-R^2
\]
we finally get equation~(\ref{first:eq}).

Since 
$\overline{AB}=2\overline{AC}\cdot\cos\gamma'_g$
and $\overline{AC}=R\cdot\cos\beta_g/\sin\gamma_g'$ we
have $\overline{AB}=2R\cdot\cos\beta_g\cot \gamma'_g$.
Substituting equations~(\ref{cotang:eq}) and~(\ref{first:eq})
in this equality we now get equation~(\ref{sec:eq}).
\finedimo

\begin{lemma}\label{diseq:lem}
We have $R'>R$.
\end{lemma}
\noindent\emph{Proof:}
By equation~(\ref{first:eq}), we have to prove that
\[
\cos^2\alpha_g+\cos^2\beta_g+2\cos^2\gamma_g+
2\cos\alpha_g\cos\beta_g\cos\gamma_g>2.
\]
Using that $\cos\gamma_g=(2\cos\alpha_g)^{-1}$
and $2\cos^2\alpha_g=1+\cos\beta_g$, this inequality 
reduces to
$2\cdot\cos^3\beta_g+5\cdot\cos^2\beta_g-1>0$,
which is verified since $\beta_g=\pi/(g+1)\geqslant
\pi/3$. 
\finedimo

Let now $O_{\Delta_g^{\!\ast}}$ be the intersection
of $\Delta_g^{\!\ast}$ with the 
horoball centered in $\infty$ and
bounded by the Euclidean plane $\matC\times\{R'\}\subset
\matC\times(0,\infty)\cong\matH^3$. By Lemma~\ref{diseq:lem}, 
$O_{\Delta_g^{\!\ast}}\subset\matC\times(0,\infty)$ 
is an Euclidean triangular prism 
and touches the three boundary triangles of $\Delta_g^{\!\ast}$.
Recall now that the geometric tetrahedra
$\Delta_g^{\! 1,\ast},\ldots,\Delta_g^{\! 2g+2,\ast}$ in $\calT_g$ 
are all isometric to $\Delta_g^{\!\ast}$. 
For any $i=1,\ldots,2g+2$ we define 
$O_i\subset\Delta_g^{\! i,\ast}$ as the only subset of
$\Delta_g^{\! i,\ast}$  
such that
the pair $(\Delta_g^{\! i,\ast},O_i)$ is isometric to
the pair $(\Delta_g^{\!\ast},O_{\Delta_g^{\!\ast}})$.
It is easily seen that the $O_i$'s glue to each other in $M_g$ 
to a horocusp section for $M_g$.
Such section is the desired maximal section
$\calO_g$, since it touches the geodesic
boundary of $M_g$.

\paragraph{Length of the slopes}
The boundary $\partial \calO_g$ of the maximal horocusp 
section of $M_g$ is a Euclidean torus tiled by $2g+2$
triangles. Let $f_1,\ldots,f_{g+1}$ be the
edges of this tiling shown in 
Fig.~\ref{tiling:fig} and 
let $\overline{s}_g$ be the slope on the boundary torus of $M_g$
determined by $f_1,\ldots,f_{g+1}$.
We set $\ell_g=L_{\calO_g}(\overline{s}_g)$.
Let us note that $\ell_g$ is the Euclidean
length of each of the $f_i$'s. 

\begin{lemma}\label{short:lem}
The following equality holds:
\[
\ell_g=2\cos\alpha_g\cdot\sqrt{\frac{4\cos^2\alpha_g-1}
{4\cos^4\alpha_g-1}}.
\]
\end{lemma}
\noindent\emph{Proof:}
Keeping notations from the previous paragraph,
we have $\ell_g=\overline{AD}/R'$. 
Now $\overline{AD}=\overline{AB}/(2\cos\gamma_g)$,
so by equation~(\ref{sec:eq}) we get 
\[
\ell_g=\frac{\cos\alpha_g+\cos\beta_g
\cos\gamma_g}{\cos\gamma_g\cdot\sqrt
{\cos^2\alpha_g+\cos^2\beta_g+\cos^2\gamma_g+2\cos\alpha_g
\cos\beta_g\cos\gamma_g-1}}.
\]
Substituting in this equation
the relations $\beta_g=2\alpha_g$, 
$\cos\gamma_g=(2\cos\alpha_g)^{-1}$ we get the desired result.
\finedimo

If $s,s'$ are slopes on a torus, we
denote by $\Delta(s,s')$ their \emph{distance},
\emph{i.e.} their geometric intersection.  

\begin{prop}\label{length:prop}
Let $g\geqslant 2$ and let $s$ be a slope on the boundary torus
of $M_g$ with $s\neq \overline{s}_g$.
Then $L_{\calO_g}(s)>6$.
\end{prop}
\noindent\emph{Proof:}
The torus $\partial\calO_g$
is tiled by $2g+2$ Euclidean triangles 
each of which has area
equal to $\ell_g^2\cdot\sin\delta_g/2$. 
Using Lemma~\ref{short:lem} and the relations
$\delta_g=\pi-2\gamma_g$, $\cos\gamma_g=(2\cos\alpha_g)^{-1}$ 
we can easily compute the area $A_g$ of $\partial\calO_g$,
getting $A_g=(2g+2)(4\cos^2\alpha_g-1)^{3/2}/(4\cos^4\alpha_g-1)$.

Let now $\theta(s,\overline{s}_g)$ be the Euclidean
angle determined by geodesic representatives for $s$ and $\overline{s}_g$
on $\partial\calO_g$. We have $L_{\calO_g}(s)\cdot
L_{\calO_g}(\overline{s}_g)\cdot\sin\theta(s,\overline{s}_g)=
\Delta(s,\overline{s}_g)\cdot A_g$, 
which implies that
\[
L_{\calO_g}(s)\geqslant \frac{\Delta(s,\overline{s}_g)\cdot A_g}
{\ell_g}=\frac{(g+1)\cdot\Delta(s,\overline{s}_g)\cdot(4\cos^2\alpha_g-1)}
{\cos\alpha_g\cdot\sqrt{4\cos^4\alpha_g-1}}.
\]
Suppose first that $(g+1)\cdot\Delta(s,\overline{s}_g)\geqslant 6$.
Then, since $4\cos^2\alpha_g=2+2\cos\beta_g\geqslant
2+2\cos \pi/3=3$, we obtain
\[
L_{\calO_g}(s)\geqslant \frac{(g+1)\cdot\Delta(s,\overline{s}_g)\cdot 2}
{\sqrt{3}}>(g+1)\cdot\Delta(s,\overline{s}_g)\geqslant 6,
\]
and we are done.

We have now to examine the cases when $\Delta(s,\overline{s}_g)=1$
and $g=2,3,4$ or $5$.
Suppose for example $g=2$, $\Delta(s,\overline{s}_2)=1$,
and let $s'_2$ be the slope corresponding to the longest
egdes of the fundamental domain for $\partial\calO_2$
shown in Fig.~\ref{tiling:fig}.
Then we have $\Delta(s'_2,\overline{s}_2)=1$,
$L_{\calO_2}(s'_2)=3\ell_2$ and $\theta(s'_2,
\overline{s}_2)=\delta_2=\arccos 1/3$. Let  
$\overline{r}_2,r'_2$ and $r$ be representatives 
in $H_1(\partial\calO_2)$ for
$\overline{s}_2,s'_2$ and $s$ respectively
(such representatives are uniquely defined up to sign).
Since $\Delta(s,\overline{s}_2)=1$, there exists
$a\in\matZ$ 
such that $r=a\cdot\overline{r}_2\pm r'_2$.
Then 
\[
L^2_{\calO_2}(s)=a^2\ell^2_2+9\ell^2_2\pm
6a\ell^2_2\cdot\cos\theta(\overline{s}_2,s'_2)=
(a^2+9\pm 2a)\cdot\ell^2_2\geqslant 8\ell^2_2. 
\]
Since $\alpha_2=\pi/6$, Lemma~\ref{short:lem} implies
that $\ell^2_2=24/5$, so $L^2_{\calO_2}(s)\geqslant
8\cdot 24/5=38.4>6^2$, and we are done.

A very similar computation works also in the remaining cases.
\finedimo

We can now prove Theorem~\ref{Dehn:teo}.
By Propositions~\ref{length:prop} and~\ref{lack:prop}
we deduce that for any $g\geqslant 2$ all but one Dehn-filling
of $M_g$ are hyperbolic. If $s^m_g$ is the meridinal slope
on the torus boundary of $M_g$, then $M_g(s^m_g)$ is the handlebody,
so it is boundary-reducible, whence non-hyperbolic. This
implies that $M_g(s)$ is hyperbolic if and only if $s\neq s^m_g$.

Moreover, if  $s\neq s^m_g$ then the Heegaard genus of $M_g(s)$
is at least $g+1$, since $\mathrm{genus}(\partial M_g(s))=g$
and $M_g(s)$ is not a handlebody. Since a Heegaard surface
for $(M_g,\partial_0 M_g,\partial_1 M_g)$ embeds in $M_g(s)$
as a Heegaard surface for $M_g(s)$ itself, the genus of
$M_g(s)$ is actually equal to $g+1$.

\section*{Acknowledgements}

The author gratefully acknowledges financial
support from the University of Melbourne and thanks Craig Hodgson 
and Damian Heard for 
many helpful conversations.

\vspace{1.5 cm}

\noindent
\hspace*{6cm}Scuola Normale Superiore\\ 
\hspace*{6cm}Piazza dei Cavalieri 7 \\
\hspace*{6cm}56127 Pisa, Italy\\ 
\hspace*{6cm}frigerio@sns.it
\vspace{.5 cm}



\end{document}